\documentclass{ifacconf}

\usepackage{graphicx}      
\usepackage{natbib}        

\usepackage{amsmath}
\newtheorem{theorem}{Theorem}

\usepackage{algpseudocode}
\usepackage{algorithm}
\usepackage{algorithmicx}

\usepackage{amsmath}
\usepackage{amssymb}
\DeclareMathOperator*{\argmin}{argmin}
 \usepackage{mathrsfs} 

\usepackage{amssymb} 
\newcommand*{\QEDB}{\hfill\ensuremath{\square}}  
\usepackage{color} 
\begin{document}
\begin{frontmatter} 

\title{Target assignment for robots constrained by limited communication 
	range\thanksref{footnoteinfo}} 

\thanks[footnoteinfo]{This work is supported in part by the China Scholarship Council (CSC).}
%
\author[Second]{Xiaoshan Bai} 
\author[First]{Weisheng Yan} 
\author[Second]{Ming Cao}
\author[Second]{Jie Huang}

\address[Second]{Faculty of Science and Engineering, University of Groningen, Groningen 9747 AG, The Netherlands \\(e-mail: \{xiaoshan.bai, m.cao, jie.huang\}@rug.nl).}
\address[First]{School of Marine Science and Technology,
	Northwestern Polytechnical University,127 West Youyi Road, Xi’an, 710072, China (e-mail: wsyan@nwpu.edu.cn).}
\begin{abstract}                
This paper investigates the task assignment problem for multiple dispersed robots constrained by limited communication range. The robots are initially randomly distributed and need to visit several target locations while minimizing the total travel time. 
A centralized rendezvous-based algorithm is proposed, under which all the robots first move towards a rendezvous position until communication paths are established between every pair of robots either directly
or through intermediate peers, and then one robot is chosen as the leader to make a centralized task assignment for the other robots. Furthermore, we propose a decentralized algorithm based on a single-traveling-salesman tour, which does not require all the robots to be connected through communication. 
We investigate the variation of the quality of the assignment solutions as the level of information sharing increases and as the communication range grows, respectively. The proposed algorithms are compared with a centralized algorithm with shared global information and a decentralized greedy algorithm respectively. 
Monte Carlo simulation results show
 the satisfying performance of the proposed algorithms. 
\end{abstract}

\begin{keyword} 
Task assignment; Robots; Limited communication range; Rendezvous-based algorithm.
\end{keyword}

\end{frontmatter}
 
\section{Introduction}
The last decade has witnessed growing demands in applying robotic networks to complete various tasks, such as terrain mapping, environmental monitoring, disaster rescue \cite{chen2010review}. The task assignment problem is on how to assign a number of robots to efficiently perform a set of tasks, which is generally managed with either centralized or decentralized algorithms \cite{korsah2013comprehensive}. 
 
Many centralized algorithms, including exact algorithms (\cite{ccetinkaya2013two}, \cite{mahmoudi2016finding}), and heuristic algorithms (\cite{edison2011integrated}, \cite{escobar2014hybrid}), were proposed to solve the task assignment problem. In \cite{dahl2009multi}, a vacancy chain scheduling was developed to formalize robot interactions for the multi-robot task assignment. Considering unmanned aerial vehicles' turning radius constraint, \cite{edison2011integrated} integrated genetic algorithms with a Dubins car model to complete target assignment for multiple aerial vehicles. Centralized algorithms can obtain optimal or near-optimal solutions for the task assignment problem; however, they require
global information \cite{choi2009consensus}. Consequently, centralized algorithms cannot solve task assignment problems in which robots only have local, or possibly outdated, information due to
the robots' limited communication capability.  

On the other hand, decentralized algorithms enable each robot to plan its own route based on available local information \cite{michael2008distributed}. For example, 
a consensus routine relying on local communication was designed for robust task allocation in \cite{choi2009consensus}. Regarding aerial vehicles with limited communication range, \cite{acevedo2014one} employed distributed approaches for surveillance mission assignment. 
In \cite{smith2009monotonic}, monotonic task assignment algorithms were proposed to minimize the time until every target is occupied by one robot with limited communication range. \cite{yu2015target} proposed decentralized algorithms to minimize the robots' total travel distance until each target position is occupied by one robot constrained by limited communication and sensing ranges. 
However, in \cite{smith2009monotonic} and \cite{yu2015target}, the numbers of robots and targets are equal, such that a robot stops moving as soon as it reaches its target.  

Motivated by the existing literature just mentioned, our research focuses on the more realistic situation, in which the number of targets is greater than the number of robots. 
To be more precise, a fleet of initially randomly distributed robots constrained by limited communication range need to visit several target locations while minimizing the total travel time. Each robot is assumed to have knowledge of all the target positions, and the positions of its communication-connected (CC) robots as in \cite{smith2009monotonic} and \cite{yu2015target}.
The main contribution of this paper is exploration of the fluctuation of the assignment quality with an increase in the quantity of information exchanged among CC robots and an increase in the robots' communication range. Firstly, we propose a centralized rendezvous-based algorithm (RBA), and a decentralized algorithm which does not require all the robots to be connected through
communication. The two algorithms enable all the target points to be visited in finite time irrespective of the robots' communication range, and the cooperative strategy used by the decentralized algorithm coordinates CC robots by effectively integrating their carried local information. Secondly, we illustrate that the quality of the solution resulted from the decentralized algorithm does not monotonically increase as the robots' communication range grows, which holds for some other decentralized task assignment algorithms. 

The rest of this paper is organized as follows. Some preliminaries are given in Section 2.  In Section 3, the formulation of the task assignment problem is presented. Section 4 studies the centralized task assignment algorithm RBA, while in Section 5 the decentralized algorithm is introduced. Monte Carlo simulations are shown in Section 6. Finally, we conclude the paper.

\section{Preliminaries}
During the robots' movement, the neighbourhoods of one robot dynamically change due to its limited communication range.
We use binary variable $c_{jk}(t)$ to model whether robot $j$ can directly communicate with robot $k$ at time $t$, namely
\begin{align} \label{eq:1}
c_{jk}(t) = \left\{
\begin{array}{lcl}
{1}, &\text{if} ~~|p_j(t)-p_k(t)|\leq r,  \\
{0}, &\text{if} ~~|p_j(t)-p_k(t)|> r,
\end{array}
\right.
\end{align}
where $p_{j}(t)$ and $p_{k}(t)$ in $\mathbb{R}^2$ are the positions of the robots $j$ and $k$ at time $t$ respectively, and $r> 0$ is the robots' communication range. 
 The following lemma shows the CC probability of a randomly distributed robotic network with respect to the limited communication range of the robots.   
\begin{lem} (Connectivity of Random Geometric Graphs, \cite{penrose1997longest})
	Consider the random geometric graph $\mathcal G$ obtained by uniformly randomly distributing $m$ points with communication range $r$ over the unit square. Then as $m\rightarrow \infty$, 
	for any $\varepsilon>0$ 
	\begin{eqnarray}\label{eq:CC}
	P(\mathcal G~ is~ connect|m\pi{r^2}-\log m\leq \varepsilon)=e^{-e^{-\varepsilon}}.
	\end{eqnarray}
	\label{Lem:1}
\end{lem}

The connectivity of the robotic communication network can be captured by the Laplacian matrix $L(t)$ of the corresponding graph as
\begin{align}
L_{jk}(t) = \left\{
\begin{array}{lcl}
{-c_{jk}(t)}, ~~~~~~~~\text{if} &j\neq k,  \\
{\sum^{m}_{i=1,i\neq j}c_{ji}(t)}, ~~\text{if} &j=k,
\end{array}
\right.
\end{align}
where $c_{jk}(t)$ is determined by (\ref{eq:1}).
 
\begin{lem} \cite{biggs1993algebraic}
	Let $\lambda_i(t), ~ i\in \{1,\cdots,m\}$, be the eigenvalues of the Laplacian matrix $L(t)$ of the multi-robot system at time $t$, where $\lambda_1(t)\leq \lambda_2(t) \leq \cdots \leq \lambda_m(t)$. Then, the multi-robot system is CC at time $t$ if and only if $\lambda_2(t)>0$.
\end{lem}

\section{Problem Statement}

We assume that in a square area with edge
length $E_l$ a set of initially randomly distributed robots, $\mathcal R=\{1,\cdots,m\}$, are employed to efficiently visit a set of dispersed target points in $\mathcal T=\{1,\cdots,n\}$. 
Each robot initially has the position information of the targets through a digital map of the environment and that of the robots within its limited communication range $r$.

The task assignment problem is to minimize total travel time of the robots accomplishing visiting all the target points. The binary variable $y_{ij}$ is applied to represent whether target $i$ is visited by robot $j$. Each robot stops moving when knowing all of its assigned targets have been visited, and starts to move again once new assignment arrives through communicating with other moving robots.
We assume that the robots move with unit speed. 
Then, the objective is to minimize 
\begin{eqnarray}\label{eq:f}
f=\sum_{j\in  \mathcal R}{t_{j}},~~~~~~~~~~~~~~~~~~~~~
\end{eqnarray}
subject to   
\begin{align}
~~~~~&\sum_{j\in \mathcal R}y_{ij}\geq1, ~~~~~~~~~~\forall~ i\in \mathcal T,~~~~~~~
\end{align}
where $t_j$ is the total travel time of robot $j$. 

\section{Centralized task assignment algorithm}  
If the robotic system is initially CC, the task assignment problem is in fact the multi-depot vehicle routing problem (MVRP) \cite{mirabi2010efficient}. The MVRP is known to be NP-hard, where a 
fleet of vehicles located at several depots need to deliver products to a
set of scattered customers. In this case, centralized algorithms are usually adopted by choosing one leader robot to make decisions for the other robots based on the global information.

Otherwise, under a certain communication range, a large number of robots can make the randomly distributed robotic system initially CC, for which Lemma \ref{Lem:1} can be applied to estimate the number of needed robots. The resulted number of robots makes the robots CC with the same probability during the whole operation time. However, in the task assignment problem, the robotic network is not necessarily always connected, since undergoing a centralized task assignment to be connected for one time is sufficient. Thus, using a large number of robots to make the robots always CC is a waste of resources.

One alternative method is to let the robots intentionally move towards a rendezvous position until the robotic system is CC. Based on this idea, we design the task assignment algorithm RBA. Inspired by the location problem of \cite{revelle2005location}, the center-of-gravity of the target points is chosen as the rendezvous position. If one robot reaches the rendezvous position first, it stops moving and waits for the other robots until all the robots are CC. 

Once the robotic system is CC, the robot having the most 1-hop neighbours is chosen to be the leader to make a centralized task assignment for all the robots, which can be dealt with using the MVRP. If there are several robots with the same number of the largest number of 1-hop neighbours, the leader is randomly chosen from these robots. In fact, each robot guided by the RBA first individually moves towards the rendezvous position, and only when all the robots are CC is a centralized task assignment made. We simply name the RBA a centralized algorithm, which is presented in Algorithm 1.

\begin{figure}[!t]                                             \centering
	\includegraphics[width=3.2 in]{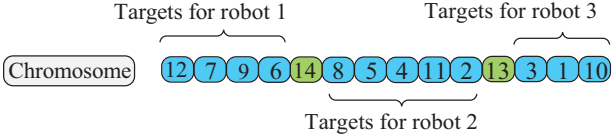}
	\caption{Chromosome structure with $12$ target points being assigned to $3$ robots, where genes $13$ and $14$ are marker genes used to separate the target points into three groups.}
	\label{fig:genestructure}
\end{figure}

The co-evolutionary multi-population genetic algorithm (CMGA) of  \cite{bai2016path} is employed for the leader robot to make the centralized task assignment. The CMGA encodes each target as a gene and inserts $m-1$ marker genes into the target genes. Thus, each chromosome represents a candidate solution to the task assignment problem. An example of the chromosome structure is presented in Fig. \ref{fig:genestructure}, which contains $12$ target points and $2$ marker genes. The routes of the $3$ robots shown in Fig. \ref{fig:genestructure} are as follows: $p_1(0)\rightarrow 12\rightarrow 7\rightarrow 9 \rightarrow 6 $, $p_2(0)\rightarrow 8\rightarrow 5\rightarrow 4 \rightarrow 11\rightarrow 2 $, $p_3(0)\rightarrow 3\rightarrow 1\rightarrow 10 $, where $p_i(0)$ is the initial position of robot $i$. As each robot is assigned with a finite number of target points, all the target points will be visited in finite time.

To study the RBA, we first introduce a lemma. Let $f_{MVA}$ be the total travel distance of the assignment solution resulted from the multi-vehicle algorithm (MVA) \cite{rathinam2007resource}, and let $f_{o}$ be the optimal value for (\ref{eq:f}).
\begin{lem} \cite{rathinam2007resource}\label{le:MVA}
	Assume that each robot initially has the position information of all the other robots, then $f_{MVA}\leq 2f_{o}$.
\end{lem} 

We use $f_{RBA}$ to represent the total travel time of the assignment solution resulted from the RBA, and $f_{MST}$ to be the sum of the distances of all edges in one
minimum spanning tree (MST) that connects all the targets and the robots.  For the pair of nodes, one being the target and the other being the robot or the two both being the targets, its edge length is the
	Euclidean distance between the two nodes, while the edge
length between nodes representing robots is $0$. 
Based on Lemma \ref{le:MVA}, we are now able to present a lower bound for the RBA when applied to solve the problem. 
\begin{theorem}\label{le:RBA}
	{\rm For the task assignment problem, the total travel time of the assignment solution resulted from the RBA is bounded by $f_{MST}\leq f_{RBA}$.	}
\end{theorem}

\begin{pf}
	The distance matrix that contains the edge length between each target point and each robot satisfies the triangle inequality. Inspired by the construction of the MVA in \cite{rathinam2007resource}, we get that Lemma \ref{le:MVA} holds as $f_{MST}\leq f_{o}$ and $f_{MVA}\leq 2f_{MST}$. 
	 As $f_{o}\leq f_{RBA}$, it holds that $f_{MST}\leq f_{RBA}$.  
 	
	 \QEDB	
\end{pf}

Theorem $\ref{le:RBA}$ can be applied to test the performance of the RBA, where a smaller $f_{RBA}/f_{MST}$ means a better performance of the algorithm. 

\alglanguage{pseudocode}
\begin{algorithm}[!tpb]
	\small 
	\caption{The centralized RBA.}
	\label{algorithm.PMA_feature}
	\begin{algorithmic}[1]
		\Require Positions of the $n$ targets, the $m$ robots and their communication range $r$.
		\Ensure Every target is visited.
		\If {The robotic network is not CC}
		\For{ Each robot $j,\forall j\in \mathcal R$}
		\If {Robot $j$ is not CC with all the other robots}		
		\State Robot $j$ moves towards the rendezvous position.
		\For {Every movement of robot $j$}
		\If {The robotic network is not CC}	
		\If{Robot $j$ does not reach the rendezvous position}
		\State Robot $j$ continues moving towards the rendezvous position.
		\Else
		\State Robot $j$ stops at the rendezvous position to wait for the other robots until the robotic network is CC. 
		\If {The robotic network is not CC}	
		\State Go to step $24$.
		\EndIf	
		\EndIf
		\Else
		\State Go to step $24$.
		\EndIf
		\EndFor
		\Else
		\State Go to step $24$.
		\EndIf
		\EndFor
		\Else
		\State One leader robot is chosen to make a centralized task assignment for the other robots, and then each robot $j$ moves to visit its assigned targets. 
		\EndIf
	\end{algorithmic}
\end{algorithm}  

\section{Decentralized Algorithm} 
If the robotic system is not initially CC, the task assignment problem can be solved in a decentralized manner where two important issues need to be considered. Those are how to plan routes for the robots which are not CC with any other robot, and how to make CC robots cooperative based on their carried partial, or even outdated, information regarding which targets have already been visited. 

As the robots are constrained by a limited communication range,  each robot $j$ carries an information tuple $\mathcal I_j(t)=\{j,p_j(t),o_j(t),t_j,u_j(t),s_j(t)\}$, where $j$ is its unique identifier; $p_j(t)$ is the current position of the robot;  $o_j(t)$, initially contains all the targets in $\mathcal T$, keeps the ordered targets on the route of robot $j$; 
$u_j(t)$ is an $m$-tuple where $u^k_j(t)=1$ means robot $j$ has ever been CC with robot $k$, otherwise $u^k_j(t)=0$; and 
$s_j(t)$ is an $n$-array used to record the target status, namely
\begin{align}
s^i_{j}(t) = \left\{
\begin{array}{lcl}
{1},~\text{if robot $j$ knows target $i$ has been visited, }~\\
{0},~\text{otherwise}.
\end{array}
\right.
\end{align}

For the task assignment problem, the worst performance of one robot happens when the robot visits all the target points without communicating with any other robot. As a result, we design a 
decentralized algorithm based on a single-traveling-salesman tour (STST), where all the targets are connected into a closed traveling salesman problem $(TSP)$ tour $TSP_0$. The $TSP_0$ is calculated by the Christofides algorithm \cite{papadimitriou1982combinatorial}, whose solution is at most 3/2 times the total travel time of an optimal TSP tour. 
Then, the robots that have never communicated with any other robot travel along the $TSP_0$ to visit the targets, which guarantees the worst performance of the robots. Thus, the first issue is solved.

The initial route $o_j(0)$ for each robot $j$ is determined by minimizing $t_j$ while traveling along the $TSP_0$, which will be shown in Theorem \ref{Theo:STSTC}.
When several robots are CC, each robot $j$ within the CC robots first updates its $n$-arrays $s_j(t)$ as
\begin{eqnarray}\label{eq:si}
s^p_j(t)=\bigcup\limits_{k\in \mathcal R_i} s^p_k(t-1),~~~\forall j\in \mathcal R_i(t), ~\forall p\in \mathcal T, 
\end{eqnarray}
where $\mathcal R_i(t)$ is the robotic subgroup $i$ in which robots are CC at time $t$. 

Then, the robots in each connected robot set $\mathcal R_i(t)$ update their remaining information tuples as follows. 
Based on the $s_j(t)$, $o_{j}(t)$ is updated by deleting the visited targets that are known by robot $j$, namely
\begin{eqnarray}
o_j(t)=o_{j}(t-1)\setminus s_{j}(t),~~~ \forall j\in \mathcal R_i(t),  
\end{eqnarray}
and then the planned total travel time $t_j$ of robot $j$ is updated to the one considering all the targets on $o_{j}(t)$ to be visited by the robot.
Afterwards, $u_j(t)$ is updated as 
\begin{align}
u^k_j(t)= \left\{
\begin{array}{lcl}
{1},~\text{if $u^k_j(t-1)=1$ or $k\in \mathcal R_i(t)$,   }~\\
{0},~\text{otherwise}.
\end{array}
\right.
\end{align}

If the robotic subgroup $\mathcal R_i(t)$ has no new members compared with $\mathcal R_i(t-1)$, the robots within the group perform their assigned tasks without task reassignment.  
To deal with how to make CC robots in $\mathcal R_i(t)$ coordinate based on their carried information, a cooperative strategy is designed. We use robotic set $\mathcal R^1_i(t)\subseteq \mathcal R_i(t)$ to contain the robots in $\mathcal R_i(t)$ that have never been CC with any other robot at time $t$, and $\mathcal R^2_i(t)$ to be $\mathcal R_i(t)\setminus \mathcal R^1_i(t)$.
 Then the targets to be divided are those in the set     
\begin{align}\label{eq:Ti}
\mathcal T_i(t) = \left\{
\begin{array}{lcl}
{o_{p}(t)},~\text{if $\sum_{k\in \mathcal R^2_i(t)}(t_k-t)\geq t_p-t$ or $\mathcal R^2_i(t)= \emptyset $, }~\\
{\bigcup\limits_{k\in \mathcal R^2_i(t)} o_k(t)},~\text{otherwise},
\end{array}
\right.
\end{align} 
where $p=\mathop{\argmin}_{q\in \mathcal R^1_i(t)}{t_q}$ if $\mathcal R^1_i(t)\neq \emptyset$, otherwise $t_p=+\infty$. 
 $\displaystyle{\sum_{k\in \mathcal R^2_i(t)}(t_k-t)\geq t_p-t}$ means that the total travel time incurred by visiting all the targets in $\bigcup\limits_{k\in \mathcal R^2_i(t)} o_k(t)$ is larger than that incurred by visiting all the targets in $o_{p}(t)$.

To make the CC robots coordinately visit the cooperative targets, the objective for the leader robot of each $\mathcal R_i(t)$ at time $t$ is to minimize
\begin{eqnarray}\label{eq:fi}
f_i(t)=\sum_{j\in \mathcal R_i(t)}{t_{j}}.
\end{eqnarray}
The resulted task assignment algorithm integrates the STST with the cooperative targets, namely (STSTC), which is shown in Algorithm $2$.

For each CC robotic network $\mathcal R_i(t)$, its leader robot employs the CMGA to make a centralized task assignment to the robots in $\mathcal R_i(t)$. If the resulted assignments do not decrease the total travel time of the robots in $\mathcal R_i(t)$,
each robot keeps its previous target assignment. 
Once a locally centralized task assignment has been completed in $\mathcal R_i(t)$, each robot $j$ updates its route $o_j(t)$.

\alglanguage{pseudocode}
\begin{algorithm}[!tpb]
	\small 
	\caption{The decentralized STSTC for each robot $j,\forall j\in \mathcal R(t)$.}
	\label{algorithm.PMA_feature}
	\begin{algorithmic}[1]
		\Require Positions of the $n$ targets, target information tuple $\mathcal I_j(t)$, communication range $r$, the tour $TSP_0$.
		\Ensure Every target on $o_j$ is visited.
		\While {$o_j(t)\neq \emptyset$ } 
		\If {Robot $j$ is not CC with any other robot}		
		\State Robot $j$ moves along its $o_j(t)$.	
		\Else
		\If {Robot $j$ is in $\mathcal R_i(t)$}
		\If { $\mathcal R_i(t)=\mathcal R_i(t-1)$}		    						
		\State Robot $j$ exchanges information with the robots in $\mathcal R_i(t)$ and continues moving without task reassignment;
		\Else
		\State Go to step $12$.
		\EndIf 	
		\EndIf
		\State One leader robot in $\mathcal R_i(t)$ is chosen to divide the cooperative targets in $\mathcal T_i(t)$ to the robots in $\mathcal R_i(t)$.
		\If {The assignments decrease the total travel time of the robots in $\mathcal R_i(t)$}
		\State Every robot $j$ in $\mathcal R_i(t)$ updates its route with $o_j(t)$.
		\Else
		\State Every robot $j$ in $\mathcal R_i(t)$ keeps its previous assignment.
		\EndIf 
	    \EndIf 
	    \EndWhile
	\end{algorithmic}
\end{algorithm} 

\subsection{Correctness of the STSTC}
To prove the correctness of the STSTC, we first present its properties. 
\begin{lem} \label{lem:STST}
	During the operation of the STSTC, the following statements hold.\\	
	~~~	(i) Each target $w\in \mathcal T$ is assigned to at least one robot, the assignment may change, but target $w$ remains assigned to at least one robot until being visited. \\	
	~~~	(ii) For robot $j$ and target $w$, if $s^w_j(t_0)=1$ at some time $t_0$, then $s^w_j(t)=1$ for all $t\geq t_0$. 
\end{lem}
 
\begin{pf} 
	Based on the initialization of the target set $o_j(t)$ for each robot $j$, the target set initially assigned to each robot is the whole target set $\mathcal T$. Thus, if one robot has never been CC with any other robot, an arbitrary target is on the robot's route unless being visited by the robot.    
	
	When several robots are CC, statement (i) is concluded based on the cooperative strategy shown in (\ref{eq:Ti}). 	
	We first consider the case when $\mathcal T_i(t)=o_p(t)$, where $p=\mathop{\argmin}_{q\in \mathcal R^1_i}{t_q}$. As robot $p$ has never been CC with any other robot, an arbitrary target $w$ satisfies $w\in o_p(t)$ if $w$ has not been visited by robot $p$. 
		
	If $\mathcal T_i(t)=\bigcup\limits_{k\in \mathcal R^2_i(t)} o_k(t)$, target $w$ is among the cooperative targets if $w\in o_k(t)$ for at least one robot $k\in \mathcal R^2_i(t)$. Otherwise, $w$ must be on the route of at least one other robot, assumed to be robot $s$, who has already communicated with at least one of the robots in $\mathcal R^2_i(t)$. As $w\in o_s(t)$, $w$ will be among the cooperative targets of a CC network if robot $s$ is CC with other robots based on the analysis when $w\in o_k(t)$ and $\mathcal T_i(t)=\bigcup\limits_{k\in \mathcal R^2_i(t)} o_k(t)$. If robot $s$ have not been CC with other robot that can visit target $w$ after wining target $w$, $w$ will be on the route of robot $s$ until being assigned to other robot or being visited.
			
	Based on the above analysis, an arbitrary unvisited target is either among the cooperative targets of a CC robotic group $\mathcal R_i$ or on the route of at least one robot. Once the target $w\in \mathcal T_i$, it will be divided among the CC robots until being visited. Otherwise, it will be assigned to at least one of the robots whose routes containing the target. Thus, statement (i) is proved.

	Statement (ii) follows directly from the union operation of the $s_j(t)$ in (\ref{eq:si}). It can also be explained by the fact that once a target is visited and its status is known by one robot, the robot will keep this information.
	
	 \QEDB	
\end{pf}

\begin{figure}[!tp] 
	\centering
	\includegraphics[height=2.2 in]{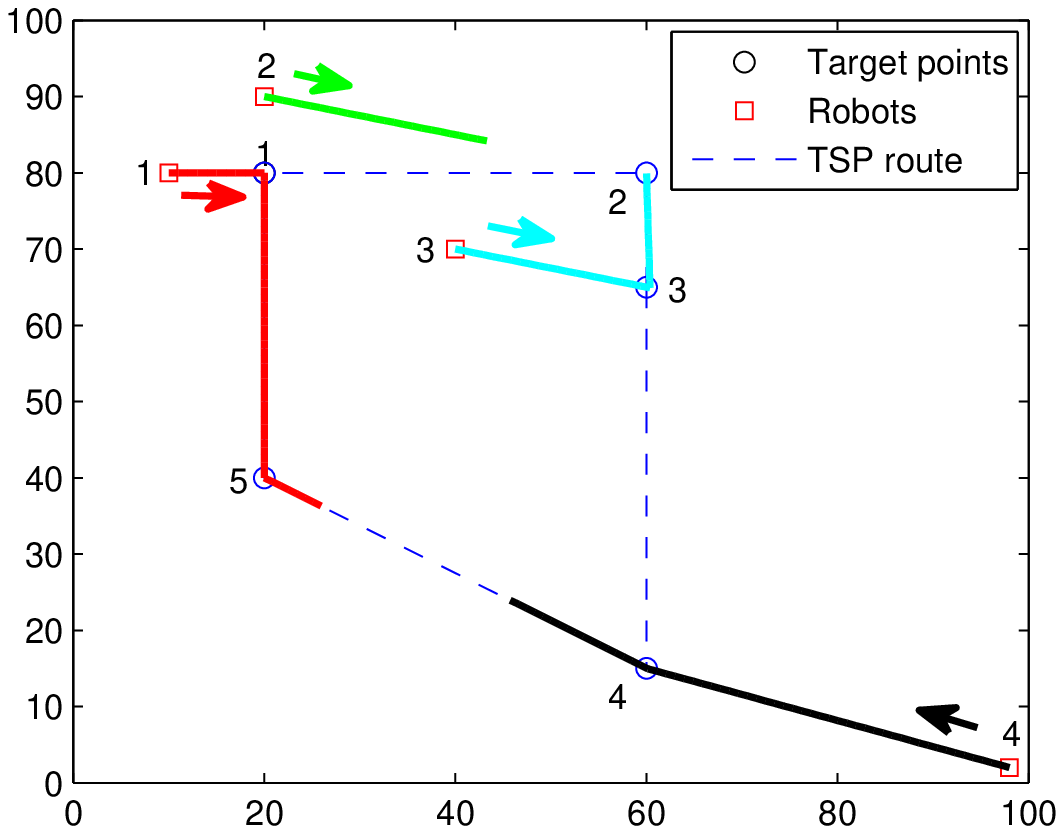}
	\caption{The routes based on the STSTC for 4 robots visiting 5 targets, where $r=25$. The total travel time is  $f=174$s.                             } 
	\label{fig:2}
\end{figure}

With these properties, we will present the main result of the section. 
\begin{theorem}
	{\rm For any given $n,m \in \mathbb{N}$, the STSTC enables all the target points to be visited in finite time. }
\end{theorem} 

\begin{pf}
	For one robot that is not initially CC with any other robot, according to the STSTC it visits targets by traveling along the tour $TSP_0$ which connects all the targets. If the robot cannot communicate with any other robot during its movement, it will not stop moving until all the targets are visited. 
	
	If the robot can communicate with other robots during its movement, the cooperative targets determined by (\ref{eq:Ti}) are divided to the CC robots. Once a target is assigned to a robot, it will be kept assigned until being visited based on Lemma \ref{lem:STST}. Thus, all the targets will be visited in finite time.  
	
	 \QEDB	
\end{pf}

We further investigate the effect of the length of the communication range of the robots on the performance of the STSTC. 
\begin{theorem}\label{theorm:2}
	{\rm For the investigated task assignment problem, a longer communication range does not necessarily lead to a better performance for the STSTC.}	
\end{theorem}

We now illustrate the main argument for this theorem. When several robots are CC, they exchange target information and update their information tuples to adjust their routes. As some targets can have already been visited by robots that are not CC, the target information shared by the CC robots is local or even outdated, which does not truly reflect the current targets' situation. As a result, the incomplete information shared by the CC robots can lead to inefficient task assignments; for example, one target that has already been visited is reassigned to another robot.  
	 	
\begin{figure}[!tp] 
		\centering
		\includegraphics[height=2.2 in]{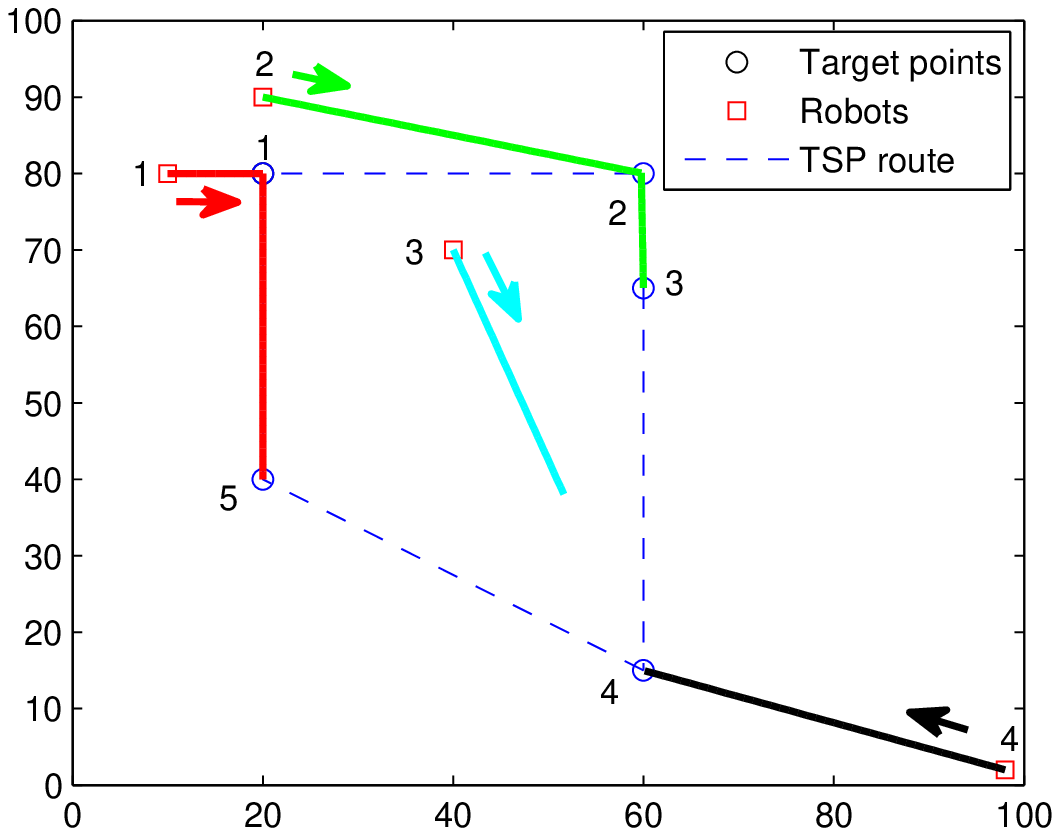}
		\caption{The routes based on the STSTC for 4 robots visiting 5 targets, where $r=30$. The total travel time is  $f=180$s.                           }
		\label{fig:3}
	\end{figure}

One illustrating case of the task assignments based on the STSTC is shown in Fig. \ref{fig:2} and Fig. \ref{fig:3} where $4$ robots with unit speed need to visit $5$ targets $\mathcal T=\{1, \cdots, 5\}$ in a $100$ $\times$ $100$ m$^2$ area. In Fig. \ref{fig:2}, the robots' communication range is $25$, which makes robots $1$ and $2$ initially CC. The routes of the $4$ robots are $p_1(0)\rightarrow 1 \rightarrow 5 \rightarrow 4$, $p_2(0) \rightarrow 2 \rightarrow 3$, $p_3(0) \rightarrow 3 \rightarrow 2 \rightarrow 1 \rightarrow 5 \rightarrow 4$, $p_4(0) \rightarrow 4 \rightarrow 5 \rightarrow 1 \rightarrow 2 \rightarrow 3$, where $p_i(0),i\in \{1,\cdots,4\}$, is the initial position of robot $i$. With the movement of the robots, target $3$ is first visited by robot $3$, and then robot $2$ stops moving once it can communicate with robot $3$ since robot $3$ is nearer to their cooperative target $2$. Reaching target $2$, robot $3$ stops moving since there is no target on $o_3$. When target $4$ is visited by robot $4$, robot $1$ should not move towards the target. However, robot $1$ continues moving as it does not have the latest target status of target $4$. When robot $1$ can communicate with robot $4$, they stop moving as their cooperative target set is empty based on (\ref{eq:Ti}). The total travel time of the robots is $174$s, where 
the travel time of the $4$ robots are $57$s, $24$s, $36$s and $57$s respectively. 
 
In Fig. \ref{fig:3}, the robots' communication range is increased to $30$, which makes robots $1$, $2$ and $3$ initially CC. After dividing the cooperative targets, we get the routes of the CC robots are $p_1(0) \rightarrow 1 \rightarrow 5$, $p_2(0) \rightarrow 2 \rightarrow 3$, $p_3(0) \rightarrow 4$, while $p_4(0) \rightarrow 4 \rightarrow 5 \rightarrow 1 \rightarrow 2 \rightarrow 3$  based on the STSTC. With the movement of the robots, robot $3$ stops moving when it can communicate with robot $4$ as robot $4$ is nearer to their cooperative target $4$. Robots $1$ and $2$ stop moving when their assigned targets are visited. The total time to visit all the targets is $180$s, where the travel time of the four robots are $50$s, $56$s, $34$s and $40$s respectively.  The performance of the STSTC for robots with a longer communication range in Fig. \ref{fig:3} is worse than that in Fig. \ref{fig:2}, thus illustrating the statement. 
Theorem \ref{theorm:2} holds for some other decentralized task assignment algorithms. 
	

\subsection{Time complexity for the STSTC}
For the robots in each $\mathcal R_i(t)$, task reassignment is centrally
made by the leader robot in the CC robotic group to minimize the total travel time of the connected robots.   
\begin{lem} \label{lem:CC}
	Each information updating for the CC robots in $\mathcal R_i$ makes the total travel time $f$ in (\ref{eq:f}) nonincreasing.   
\end{lem}

\begin{pf}
	For the STSTC, the CC robots in each $\mathcal R_i$ will divide the cooperative targets to the robots if the reassignment reduces their total travel time (\ref{eq:fi}). Otherwise, each robot $j$ among the CC robots in each $\mathcal R_i$ only updates its $s_j$, $o_j$, $o_j$, and the corresponding $t_j$ based on their shared information, which does not make the total travel time $f_i$ in (\ref{eq:fi}) worse. As $f_i$ is one component of $f$ in (\ref{eq:f}), $f$ is nonincreasing for each information updating. 
	
	\QEDB	
\end{pf}

For the robot that is not initially CC with any other robot, it will visit all the targets through circling the $TSP_0$ if it cannot communicate with any other robot during its movement. The following lemma gives an upper bound on the minimal travel distance of the robot. 
\begin{lem} \cite{supowit1983travelling}\label{lem:1}
	Let $L$ be the shortest length of the TSP tour which optimally connects the $n$ target points and one robot in a square area with edge length $E_l$, then there exist $c_1\in \mathbb{N}$ and $c_2\in \mathbb{R}^{+}$ such that $L\leq c_2\sqrt{(n+1)*E_l}$ for all $n+1\geq c_1$.
\end{lem}
As the robots move with unit speed, the upper bound on the minimal travel time for one robot visiting the $n$ targets is not larger than $c_2\sqrt{(n+1)*E_l}$.
With this lemma, we are able to give an upper bound on the time complexity of the STSTC, where time complexity means that all the target points are visited in finite time.   
\begin{theorem}\label{Theo:STSTC}
	{\rm The STSTC guarantees the $n$ target points to be visited by $m$ robots with the total travel time at most $m(\sqrt{2E_l}+(3/2)c_2\sqrt{nE_l})$ where $c_2$ is obtained in Lemma \ref{lem:1}. }
\end{theorem} 

\begin{pf} 
	We use matrix $D$ as the distance matrix where $D(i,j)$ contains the distance between vertex $i$ and $j$. $L(\{p_i(0),T^{1}_i, T^{2}_{i}, \cdots, T^{n}_i\})$ is employed to represent the length of the initial route $o_i(0)$ for robot $i$, where $o_i(0)$ is $p_i(0) \rightarrow T^{1}_i \rightarrow T^{2}_{i} \rightarrow \cdots \rightarrow T^{n}_i$ based on the $TSP_0$. Then
	\begin{eqnarray}
	L(o_i(0)) &=& D(p_i(0),T^{i}_i)+D(T^{1}_i,T^{2}_i)+ \cdots  \nonumber\\
	&&+D(T^{n-1}_{i},T^{n}_i) \nonumber\\
	&=& D(p_i(0),T^{1}_i)+D(T^{1}_i,T^{2}_i)+ \cdots \nonumber \\
	&&+ D(T^{n-1}_{i},T^{n}_i) 
	+D(T^{n}_{i},T^{1}_i)-D(T^{n}_{i},T^{1}_i) \nonumber \\
	&=& D(p_i(0),T^{1}_i)+L(\{T^{1}_i, \cdots, T^{n}_i,T^{1}_i\}) \nonumber \\
	&&     -D(T^{n}_{i},T^{1}_i) \nonumber \\
	&= & D(p_i(0),T^{1}_i)+L(TSP_0)-D(T^{n}_{i},T^{1}_i),\nonumber \\
	\end{eqnarray}
	where $L(TSP_0)$ is the length of the tour $TSP_0$.   		
	As $D(p_i(0),T^{1}_i) \leq \sqrt{2E_l}$ and $L(TSP_0)\leq (3/2)c_2\sqrt{nE_l}$ where $TSP_0$ is calculated by the Christofides algorithm \cite{papadimitriou1982combinatorial}, we get $L(o_i(0))\leq \sqrt{2E_l}+(3/2)c_2\sqrt{nE_l}$. If robot $i$ cannot communicate with any other robot during its movement, it will travel along the  $o_i(0)$ until all the targets being visited.
	Thus, the longest travel distance of one robot guided by the STSTC is upper bounded by $\sqrt{2E_l}+(3/2)c_2\sqrt{nE_l}$. Moreover, the first target to be visited, $T^{1}_i$, and the travel direction of the robot $i$ when traveling along the $TSP_0$ are chosen based on the minimization of $L(o_i(0))$. In other words, $o_i(0)$ is determined by minimizing $L(o_i(0))$.

	Based on Lemma \ref{lem:CC}, the worst performance of the STSTC occurs when each robot circles around the tour $TSP_0$ without communicating with any other robot.   
	In this case, each robot stops moving after visiting all the targets on $TSP_0$ by itself. Since all the robots move with normalized unit speed, the total travel time is the total travel length of the robots, thus proving the theorem. 
	
	 \QEDB	
\end{pf}

\section{Monte Carlo study}

We implemented two sets of simulations in a $1000$ $\times$ $1000$ m$^2$ area where the numbers of target points and robots are $15$ and $4$, $30$ and $6$ respectively. For each set of simulations, Monte Carlo simulations are carried out on $500$ scenarios where the positions of the target points and the robots are randomly generated. 

The proposed algorithms are compared with a greedy algorithm where robots always move towards the nearest target and the CC robots communicate only when two or more robots are moving towards the same target. The assigned results are also compared with the centralized Christofides algorithm \cite{papadimitriou1982combinatorial}, where a lower bound on the solution of the task assignment problem is obtained based on the global information. 
All the experiments are performed on an Intel Core (TM) $i5-$ $4590$ CPU $3.30$ GHz, with algorithms compiled by Matlab under Windows $7$.  

The solution quality of each algorithm is defined by  
\begin{eqnarray}\label{eq:Compqare}
q&=&\frac {f}{f_{MST}},
\end{eqnarray}
where $f$ is the value in (\ref{eq:f}); $f_{MST}$ is a tight lower bound of the solution calculated based on the MST \cite{papadimitriou1982combinatorial} where all the global information is assumed to be available. Thus, a smaller $q$ of one algorithm means a better performance of the algorithm. 
\begin{figure}[!tp] 
	\centering
	\includegraphics[height=2.2 in]{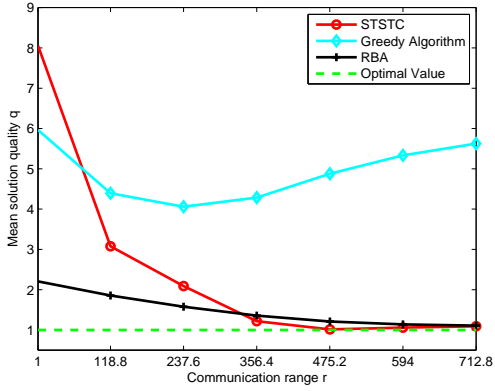}
	\caption{The mean solution quality for $4$ robots with different communication ranges visiting $15$ targets, where $P(connect~|~r_c=712.8m)=0.99$.                             } \label{fig:N15k3}
\end{figure}

\begin{figure}[!tp] 
	\centering
	\includegraphics[height=2.2 in]{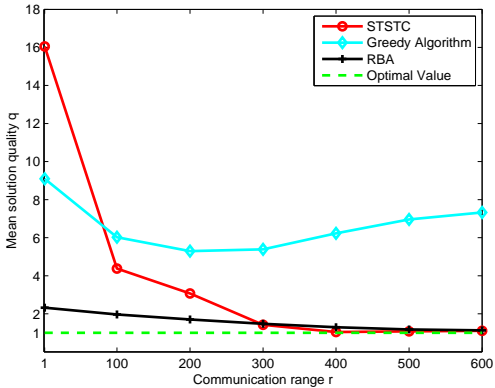}
	\caption{The mean solution quality for $6$ robots with different communication ranges visiting $30$ targets, where $P(connect~|~r_c=600m)=0.99$.  (change the graph!)                     } 
	\label{fig:N30k6}
\end{figure}

Testing the algorithms on robots constrained by communication range varying from $1$m to the one that makes the whole robotic system CC with probability $0.99$, we show the assignment results of the two sets of simulations in Fig. \ref{fig:N15k3} and Fig. \ref{fig:N30k6}. When the robots' communication range $r$ is $0$, information sharing exists among robots whose current positions are the same. The two figures first show that the quality of the assignment solutions resulted from the decentralized STSTC and the greedy algorithm varies greatly as $r$ grows, where a
longer $r$ generally leads to a better performance for the STSTC but not for the the greedy algorithm. The reason is that more environmental information is shared among the CC robots guided by the STSTC where cooperative targets are properly divided, while the CC robots guided by the greedy algorithm communicate only when task assignment conflicts occur. However, the greedy algorithm outperforms the STSTC when $r$ is approximately smaller than $r_c/6$ where $P(connect~|~r_c)=0.99$. That is because the robots using the STSTC cannot cooperate with other robots frequently when $r$ is short, while the robots guided by the greedy algorithm do not rely on the communication so much and have relatively smaller travel cost by always moving towards the nearest target.

Furthermore, Fig. \ref{fig:N15k3} and Fig. \ref{fig:N30k6} show that the RBA has a better performance which does not vary greatly with an increase in $r$. The reason is that the longer $r$ only makes all the robots guided by the RBA CC at a earlier time, and then they cooperatively visit the remaining targets. Though the STSTC does not perform well when $r$ is short, it outperforms the RBA when $r\geq r_c/2$, as shown in the two figures. As
for the STSTC, a longer $r$ leads to more information shared by the CC robots, which generally results in better cooperation for the CC robots. However, a worse performance of the the STSTC occurs in Fig. \ref{fig:N30k6} when $r$ is increased from $400$ to $500$, which is generally the case for the decentralized greedy algorithm as shown in the two figures. Thus, Theorem \ref{theorm:2} is again verified.
  
  \begin{figure}[!tp] 
  	\centering
  	\includegraphics[height=2.2 in]{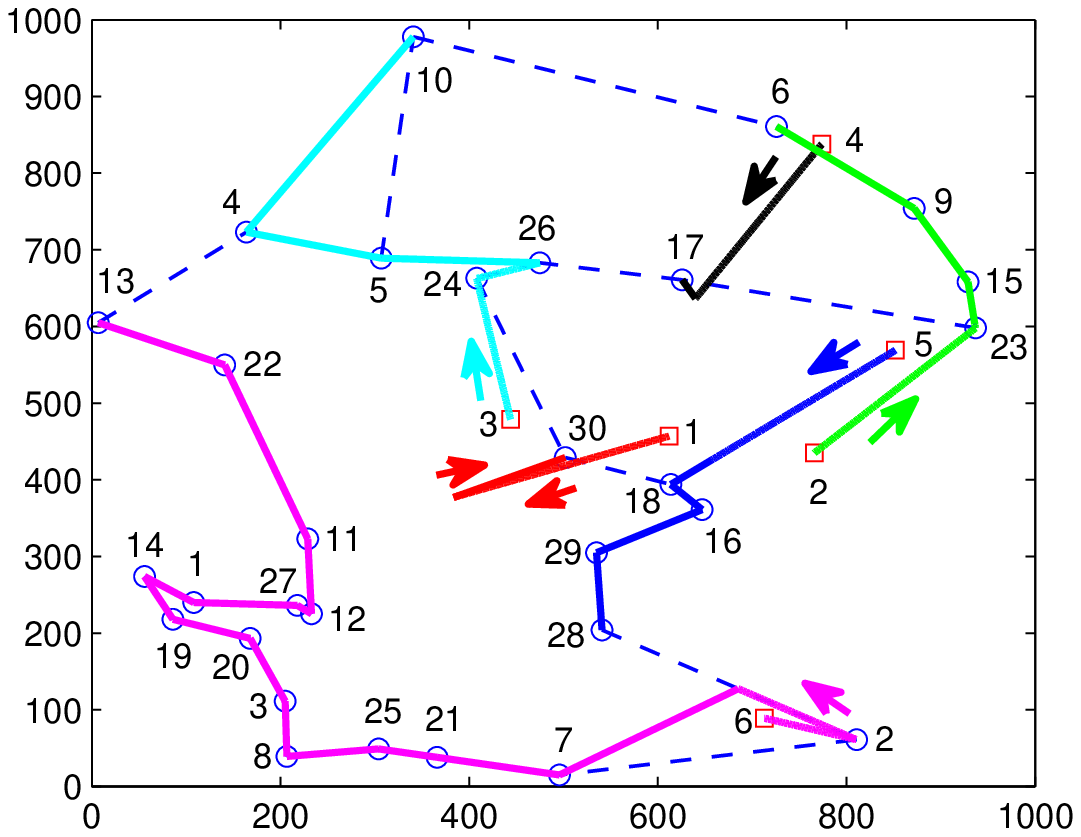}
  	\caption{The routes based on the STSTC for $6$ robots visiting $30$ target points, where $r=300$m. The total travel time is $2973.1$s.                             } \label{fig:STSTAN30k6}
  \end{figure}
  
  \begin{figure}[!tp]  
  	\centering
  	\includegraphics[height=2.2 in]{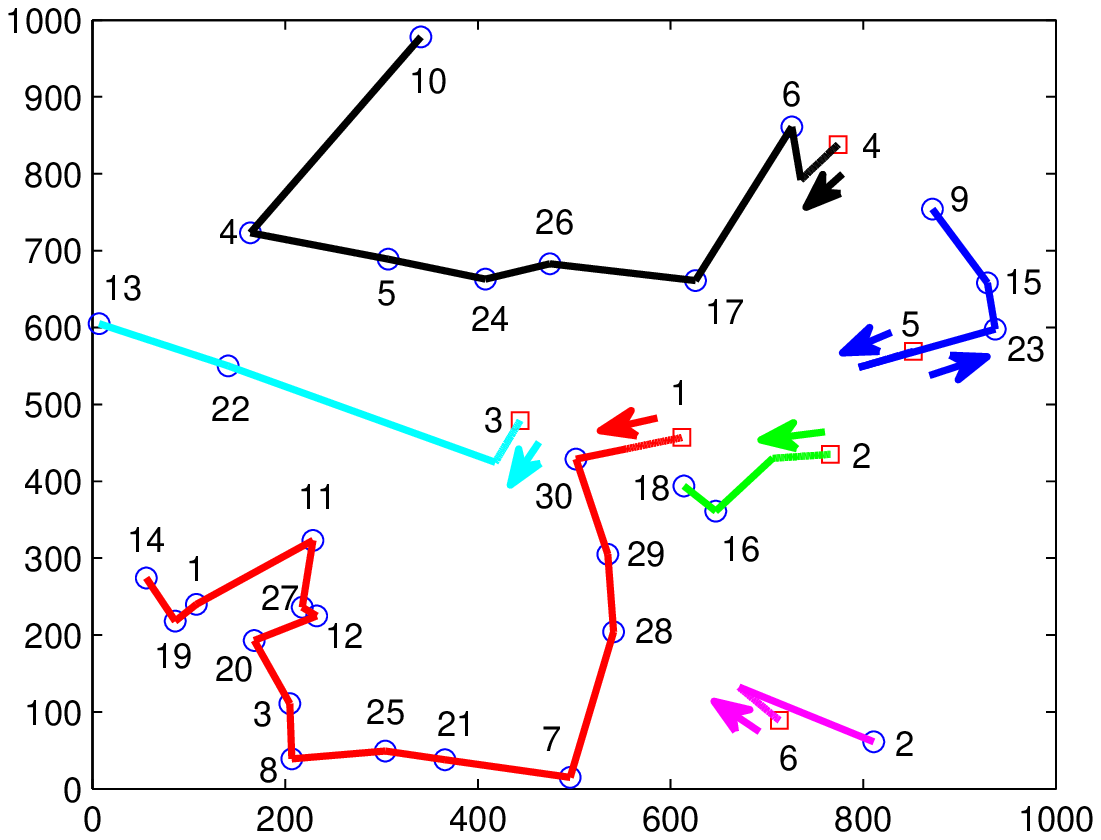}
  	\caption{The routes based on the RBA for $6$ robots visiting $30$ target points, where $r=300$m. The total travel time is $3854.5$s.                            } 
  	\label{fig:RBAN30k6}
  \end{figure}
  
Finally, the mean solution quality $q$ of the RBA and the STSTC displayed in Fig. \ref{fig:N15k3} is approximate to the optimal value $1$ when $r\geq r_c/2$, which shows the good performance of the algorithms. 
However, in Fig. \ref{fig:N30k6}, $q$ is a little bit larger
than $1$ when $r\geq r_c/2$. There are two reasons: one being that the $f_{MST}$ in (\ref{eq:Compqare}) is a tight lower bound of the optimal value, which can lead to a larger $q$; the other one being that the calculation of the $f_{MST}$ is based on the global information, while the assignments resulted from the RBA and the STSTC consider the limited communication range $r$ among the robots. 

We also show the routes, resulted from the STSTC and RBA, for $6$ robots visiting $30$ target points in one scenario in Fig. \ref{fig:STSTAN30k6} and Fig. \ref{fig:RBAN30k6}, where $r=r_c/2$. In Fig. \ref{fig:STSTAN30k6}, the routes of the robots generally follow a single TSP tour which connects all the targets. The robots in Fig. \ref{fig:RBAN30k6} are guided by the RBA in which all them first move towards the center-of-gravity of the targets to make them CC. Then, every robot moves to their assigned target points after a centralized task assignment is made. As the robots need to intentionally move towards the rendezvous position to make them CC, the performance of the RBA in the experimental scenarios is worse than that of the STSTC. The result verifies the effectiveness of the STSTC when $r\geq r_c/2$ as shown in Fig. \ref{fig:N15k3} and Fig. \ref{fig:N30k6}.  

It can be concluded from the above analysis that the centralized RBA leads to satisfactory assignment solutions irrespective of the robots' communication range $r$, while the decentralized STSTC has competitive performance when $r\geq r_c/2$. However, relying on local information makes the STSTC more robust than the centralized RBA if robot failure is concerned. 

\section{Conclusion}
In this paper, we have studied the task assignment problem where several initially randomly distributed robots constrained by limited communication range are coordinated to visit a set of dispersed target points. The centralized algorithm RBA and the decentralized STSTC proposed in the paper guarantee all the target points to be visited in finite time irrespective of the communication range. For the task assignment problem, we have illustrated that a longer communication range does not necessarily lead to a better performance for the STSTC, which usually holds for the other decentralized algorithms, while Monte Carlo simulations have shown that longer communication ranges lead to better performances for the RBA. The proposed algorithms will be further tested by considering environmental disturbances, for example, winds and obstacles. We are also planning to test the algorithms using real mobile robots.
  

\bibliography{wenxian}                                                     

\end{document}